\documentclass[11pt,twoside]{article} 
\usepackage{latexsym} 
\usepackage{enumitem}
\usepackage{color}

\input{amssym.def} 
\input{amssym} 
\newfont{\fra}{eufm10 scaled 1095} 
\newfont{\Bb}{msbm10 scaled 1095} 
\newfont{\Bbg}{msbm10 scaled 1680} 
\newcommand\CC{{\mbox{\Bb C}}} 
\newcommand\RR{{\mbox{\Bb R}}} 
\newcommand\NN{{\mbox{\Bb N}}} 
\newcommand\ZZ{{\mbox{\Bb Z}}} 
\newcommand\QQ{{\mbox{\Bb Q}}}

\newcommand\fg{{\frak{g}}} 
\newcommand\fh{{\frak h}}

\newcommand\fl{{\frak l}} 
 
\newcommand\fn{{\frak n}}

\newcommand\fa{{\frak a}} 
\newcommand\fd{{\frak d}}

\newcommand\fz{{\frak z}} 
\newcommand\dt{{\hat t\hspace{1pt}}}

\newcommand\cA{{\cal A}}

\newcommand{\fsl}{\mathop{{\frak s \frak l}}} 
\newcommand\fosc{\frak o \frak s \frak c} 
 
\newcommand{\fso}{\mathop{{\frak s \frak o}}}

\newcommand{\GSp}{\mathop{{\rm GSp}}} 
\newcommand{\SL}{\mathop{{\rm SL}}} 
\newcommand{\GL}{\mathop{{\rm GL}}} 
\newcommand{\grO}{{{\rm O}}}
\newcommand{\Sp}{\mathop{{\rm Sp}}} 
\newcommand{\Hom}{\mathop{{\rm Hom}}} 
\newcommand{\Aut}{\mathop{{\rm Aut}}} 
\newcommand{\Osc}{{{\rm Osc}}} 
\newcommand{\dH}{{\rm H}_{q+1}(\ZZ)} 
\newcommand{\Id}{{{\rm id}}} 
\newcommand{\ad}{{{\rm ad}}}

\newcommand{\Ad}{\mathop{{\rm Ad}}}

\newcommand{\Ker}{\mathop{{\rm ker}}}

\renewcommand{\Re}{\mathop{{\rm Re}}} 
\renewcommand{\Im}{\mathop{{\rm im}}} 
\newcommand{\diag}{\mathop{{\rm diag}}}

\newcommand{\Span}{{{\rm span}}} 
\newcommand{\mod}{\mathop{{\rm mod}}} 
\newcommand{\proj}{{{\rm pr}}} 
\newcommand\ip{{\langle\cdot \,,\cdot \rangle}}

\newcommand\la{{\langle}}
\newcommand\ra{{\rangle}}
\newcommand\proof{{\sl Proof. }} 
\newcommand{\qed}{\hspace*{\fill}\hbox{$\Box$}\vspace{2ex}} 
\newcommand{\benur}{\begin{enumerate}[label=(\roman*)]}
\newtheorem{theo}{Theorem}[section] 
\newtheorem{pr}[theo]{Proposition}

\newtheorem{ex}[theo]{Example}
\newtheorem{re}[theo]{Remark}
\newtheorem{co}[theo]{Corollary}
\newtheorem{lm}[theo]{Lemma}
\begin{document} 
\title{Existence of cocompact lattices in Lie groups with a bi-invariant metric of index 2} 
\author{Ines Kath}
\maketitle 
\begin{abstract}
\noindent
We study the existence of cocompact lattices in Lie groups with bi-invariant metric of signature $(2,n-2)$. We assume in addition that the Lie groups under consideration are simply-connected, indecomposable and solvable. Then their centre is one- or two-dimensional. In both cases, a parametrisation of the set of such Lie groups is known. We give a necessary and sufficient condition for the existence of a lattice in terms of these parameters.
For groups with one-dimensional centre this problem is related to Salem numbers.
\end{abstract}

\section{Introduction} 

The main motivation for this paper is the following result by Baues and Globke \cite{BG}: Let $M$ be a compact pseudo-Riemannian manifold, and let $G$ be a connected solvable Lie group of isometries acting transitively on $M$. Then the stabiliser of any point $x \in M$ is a discrete subgroup in $G$. In particular, $G\to M,\ g\mapsto g\cdot x$ is a covering map and $M$ is diffeomorphic to a quotient of $G$ by a cocompact lattice. The pseudo-Riemannian metric g on $M$ pulls back to a left-invariant metric $g_G$ on $G$. Baues and Globke proved that $g_G$ is actually bi-invariant. Hence, if one is interested in compact pseudo-Riemannian homogeneous spaces of solvable Lie groups, then one is lead to the question which Lie groups with bi-invariant metric have cocompact lattices. Obviously, one may assume that $G$ is simply-connected. Furthermore, we want to assume that $(G,g_G)$ is indecomposable, although, in general, a lattice in a direct product $G_1\times G_2$ of Lie groups is not a product of lattices in $G_1$ and $G_2$. 

There is a one-to-one correspondence between simply-connected Lie groups with a bi-invariant metric and metric Lie algebras (also called quadratic Lie algebras). 
By a metric, we always mean a non-degenerate but not necessarily definite inner product. The index of a metric is defined as the number of minus signs in its Sylvester normal form. 
By the index of a metric Lie algebra $(\fg,\ip)$ we mean the index of the invariant inner product $\ip$. For index one or two, every indecomposable metric Lie algebra either is simple or solvable. Borel proved that each simple Lie group has a cocompact lattice \cite{B}. Besides our main motivation mentioned above this is a further reason why we will concentrate here on lattices in solvable Lie groups. 

Solvable Lorentzian metric Lie algebras were classified by Medina \cite{M85}. It turns out that these Lie algebras are exactly the oscillator algebras. In \cite{MR85}, Medina and Revoy answer the question which of the associated Lie groups admit a lattice. However, notice that the classification of such lattices in \cite{MR85} is far from being correct. A description of lattices in oscillator groups and a classification in dimension four was achieved by Fischer \cite{F}.

Here we consider the case of index two. A classification of metric Lie algebras of index two was given in \cite{BK}. It was redone in \cite{KO0}, where we used a much more systematic approach. The centre $\fz(\fg)$ of an indecomposable solvable metric Lie algebra $(\fg,\ip)$ of index two is one- or two-dimensional. If $\dim\fz(\fg)=1$, then $\fg$ is even-dimensional and isomorphic to a generalised oscillator algebra. The set of such algebras is parametrised by $q=(\dim \fg -2)/2$ real numbers $\mu_1,\dots,\mu_q\in\RR\setminus\{0\}$. If $\dim \fz(\fg)=2$ and if $\dim \fg$ is even, then the structure of $\fg$ is similar to that of an oscillator algebra. If $\dim \fg$ is odd, then the structure of $\fg$ is slightly more complicated.  In both cases the set of such Lie algebras is parametrised by $q=[(\dim \fg -4)/2]$ linear forms $\mu_1,\dots,\mu_q\in (\RR^2)^*\setminus\{0\}$.

For each of these metric Lie algebras of index two we describe the associated simply-connected Lie group and we give a 
necessary and sufficient
criterion in terms of $\mu_1,\dots,\mu_q$ for the existence of a lattice in this group. If the centre of $G$ is one-dimensional, this criterion uses Salem numbers. If $f$ is the minimal polynomial over $\QQ$ of a Salem number $r>1$, then the roots of $f$ different from $r^{\pm1}$ lie on the unit circle. Let $Z_r$ denote the multiset of these roots. We will see that, essentially,
$G$ admits a lattice  if and only if there exists a Salem number $r>1$ such that $Z_r$ consists of all numbers $r^{ \pm i \mu_1},\dots,r^{ \pm i \mu_q}$ different from 1. For a precise formulation see Theorem~\ref{T1}. If the centre is two-dimensional, then there exists a lattice in $G$ if and only if the linear forms $\mu_1,\dots,\mu_q$ are contained in a lattice of $(\RR^2)^*$, see Theorem~\ref{T2}.

\section{Lie groups with bi-invariant metric of index 2}\label{S2}
\subsection{Basic facts}
Let $G$ be a Lie group and $g$ be a semi-Riemannian bi-invariant metric on $G$.
The infinitesimal object associated to $(G,g)$ is a pair $(\fg,\ip)$ consisting of the Lie algebra $\fg:=T_eG$ of $G$ and the $\fg$-invariant inner product $\ip:=g_e$. Any such pair $(\fg,\ip)$ is called a {\it metric Lie algebra}. Some authors use the term {\it quadratic Lie algebra} instead of metric Lie algebra.
We will say that two Lie groups with bi-invariant metric $(G_1, g_1)$ and $(G_2,g_2)$ are {\it isomorphic} if there is a group isomorphism $\phi:G_1\rightarrow G_2$ that is also an isometry.  Isomorphisms of metric Lie algebras are defined analogously. A Lie group with bi-invariant metric is called {\it indecomposable} if it is not isomorphic to a direct product $(G_1\times G_2, g_1\times g_2)$ of at least one-dimensional Lie groups with bi-invariant metric. Analogously, a metric Lie algebra is called {\it indecomposable} if it is not the direct sum of non-trivial metric Lie algebras. 

In \cite{KO0, KO1} a classification scheme for metric Lie algebras is discussed. Every metric Lie algebra whose underlying Lie algebra contains a semisimple ideal either is semisimple or is the direct sum of a semisimple metric Lie algebra and one without simple ideals. In \cite{KO1} it is shown that a metric Lie algebra without simple ideals has the structure of a quadratic extension. Moreover, a description of  quadratic extensions in terms of a quadratic cohomology is given. For small index of the metric, this cohomology set can be calculated explicitly, which gives an explicit classification in this case. For index $2$, a classification was already obtained in \cite{KO0}, where a special case of quadratic extensions is discussed.

\subsection{Metric Lie algebras of index 2} \label{Sqe}

Let $(\fg, \ip)$ be an indecomposable metric Lie algebra of signature $(2, q)$. If $\fg$ is simple, then $\fg$ is isomorphic to $\fsl(2,\RR)$ and $\ip$ is a multiple of the Killing form. Since it is well known that $\SL(2,\RR)$ has cocompact lattices, we will concentrate on non-simple indecomposable metric Lie algebras
of index 2. These Lie algebras are solvable.

As explained above, all  these Lie algebras have the structure of a quadratic extension. Here we will not need the notion of a general quadratic extension. All relevant  metric Lie algebras in this paper are quadratic extension of the following special kind. We take an abelian Lie algebra
$\fl\in\{\RR,\RR^2\}$ and its dual $\fz:=\fl^*$. Let $(\rho,\fa,\ip_\fa)$ be an orthogonal $\fl$-module, 
i.e. a triple consisting of a (pseudo-)Euclidean vector space $(\fa,\ip_\fa)$ and an orthogonal representation $\rho:\fl \to \fso(\fa,\ip_\fa)$,
and let  $\alpha\in \Hom(\bigwedge^2\fl,\fa)$ be a 2-cocycle.
We endow the vector space $\fg:=\fz\oplus\fa\oplus\fl$ with a non-degenerate inner product $\ip$ given by
$$
\la (z,a,t), (z',a',t')\ra:= z(t')+z'(t)+\la a, a'\ra_\fa
$$
and a Lie bracket such that $[\fz,\fg]=0$ and 
\begin{eqnarray}
\,[l_1,l_2]&=&\alpha(l_1,l_2)\in\fa,  \label{c1}\\
\, [l,a]&=& -\la a, \alpha(l,\cdot)\ra+\rho(l)(a) \in\fz\oplus\fa,\label{c2}\\
\, [a_1,a_2]&=&\la \rho(\cdot) a_1,a_2\ra \in\fz\nonumber
\end{eqnarray}
for all  $l,l_1,l_2\in\fl$ and $a,a_1,a_2\in \fa$. Then $(\fg,\ip)$ is a metric Lie algebra.

We consider the following special choices of $\fl$, $\fa$, $\rho$ and $\alpha$:

\benur
\item
Consider $\fl=\RR$ and let $(\fa,\ip_\fa):=\RR^{1,1+2q}=\RR^{1,1}\oplus \RR^{2q}$ be the $(2q+2)$-dimensional Minkowski space, $q\ge0$.
We fix an element $\mu=(\mu_1,\dots, \mu_q)\in (\RR\setminus\{0\})^q$ if $q\ge1$. For $q=0$, we put $\mu=\emptyset$. Identifying $\fa=\RR^{1,1+2q}\cong \CC^{1+q}$ as vector spaces, we define a linear map 
$L\in \fso(\fa,\ip_\fa)$
by 
$$L:\CC^{1+q}\longrightarrow \CC^{1+q},\quad (z_0,z_1,\dots,z_q)\mapsto i(\overline{z_0},\mu_1 z_1,\dots, \mu_q z_q)$$
and an orthogonal representation $\rho$ of $\RR$ on $\fa$ by $\rho(t):=t L$. Put $\alpha=0$.
The arising metric Lie algebra is a generalised oscillator algebra and will be denoted by $\fosc_{1,q}(1,\mu)$. 
\item
Now take $\fl=\RR^2$ and let $(\fa,\ip_\fa)=:\fa^{2q}$ be the $2q$-dimensional Euclidean space. We fix elements $\mu_1,\dots,\mu_q$ in the dual space of $\fl$ and define an orthogonal representation $\rho$ of $\fl$ on $\fa^{2q}\cong \CC^q$ by
$$\rho(t)(z_1,\dots, z_q)= i(\mu_1(t) z_1,\dots, \mu_q(t) z_q)$$
for $t\in\RR^2$. Put $\alpha=0$.
The arising metric Lie algebra will be denoted by
$\fosc^{2}_{q}(\mu)$.

\item Let $\fl$, $\fa^{2q}$ and $\rho$  be defined as in (ii).  We put $\fa=\fa^{2q}\oplus \fa_0$, where $\fa_0=\RR^1$ endowed with the standard metric. We consider $\rho$ as a representation on $\fa$ being trivial on $\fa_0$. We define a 2-form $\alpha: \bigwedge^2\fl \rightarrow \fa_0$ by $\alpha(e_1,e_2)=1$, where $e_1,e_2$ is the standard basis of $\fl=\RR^2$.  We denote the arising metric Lie algebra by
$\fd_q(\mu)$.

\end{enumerate}
The symmetric group on $q$ letters, $S_q$, acts on $\RR^q$ by permuting the coordinates. Moreover, $(\ZZ_2)^q$ acts on $\RR^q$ by component-wise multiplication. This action extends to an action on subspaces of $\RR^q$ and on alternating forms on $\RR^q$. Let $e_1,e_2$ be the standard basis of $\RR^2$ and put $\mu(\RR^2):=\Span\{\mu(e_1),\mu(e_2)\}$ for $\mu\in ((\RR^2)^*)^q$. 
\begin{pr}\label{P1}{\rm \cite{KO0}} Let $(\fg,\ip)$ be a non-simple indecomposable metric Lie algebra of signature $(2,n-2)$. Then the centre $\fz(\fg)$ of $\fg$ is one- or two-dimensional and we are in one of the following cases. 
\begin{enumerate}
\item If $\dim \fz(\fg)=1$, then $n$ is even and $\fg$ is isomorphic to $\fosc_{1,q}(1,\mu)$  for exactly one $\mu=(\mu_1,\dots, \mu_q)\in \RR^q$ with $0<\mu_1\le\dots \le\mu_q$, where $q:=n/2-2$.
\item If $\dim \fz(\fg)=2$ and $n$ is even, then $q:=n/2-2\ge 3$ and $\fg$ is isomorphic to $\fosc^2_{q}(\mu)$ for some $\mu=(\mu_1,\dots,\mu_q)\in ((\RR^2)^*\setminus 0)^q$.  The set $\{\mu_1,\dots,\mu_q\}\subset (\RR^2)^*$ is not contained in the 
    union of two 1-dimensional subspaces.
    Furthermore, $\fosc^2_{q}(\mu)$ and $\fosc^2_{q}(\mu')$ are isomorphic as metric Lie algebras if and
    only if
    $$\mu(\RR^2)=\mu'(\RR^2)\quad \mod  
    S_{q}\ltimes(\ZZ_{2})^{q}.$$
\item If $\dim \fz(\fg)=2$ and $n$ is odd, then $q:= (n-5)/2\ge 0$ and $\fg$ is isomorphic to $\fd_{q}(\mu)$ for some $\mu=(\mu_1,\dots,\mu_q)\in((\RR^2)^*\setminus 0)^q$.
    Furthermore, $\fd_{q}(\mu)$ and $\fd_{q}(\mu')$ are isomorphic as metric Lie algebras if and
    only if
$$(\,\mu(\RR^2),\ \pm\mu(e_1)\wedge \mu(e_2)\,) =(\,\mu'(\RR^2),\ \pm\mu'(e_1)\wedge \mu'(e_2)\,)\quad \mod  S_{q}\ltimes(\ZZ_{2})^{q}.$$
\end{enumerate}
Conversely, if $(\fg,\ip)$ is one of the metric Lie algebras listed in items 1., 2., or 3., then $\ip$ has index 2 and $(\fg,\ip)$ is indecomposable. 
\end{pr}

The preceding proposition immediately implies a classification of indecomposable simply-connected Lie groups with bi-invariant metric of signature $(2,n-2)$. Since the 
existence of lattices is independent of the metric we are interested in isomorphism classes of the underlying Lie groups. Let us discuss this problem on the level of Lie algebras.
\begin{pr} \label{P2}
The Lie algebras $\fosc_{1,q}(1,\mu)$ and $\fosc_{1,q}(1,\mu')$ are isomorphic if and only if they are isomorphic as metric Lie algebras. The same is true for $\fosc^2_{q}(\mu)$ and $\fosc^2_{q}(\mu')$. The Lie algebras  $\fd_{q}(\mu)$ and $\fd_{q}(\mu')$ are isomorphic if and only if  $$\mu(\RR^2)=\mu'(\RR^2)\ \mod  
    S_{q}\ltimes(\ZZ_{2})^{q}.$$
\end{pr}
\proof
Consider first the Lie algebra $\fg:=\fosc_{1,q}(1,\mu)$. Then the Lie algebras $[\fg,\fg]/\fz$ and $\fg/[\fg,\fg]$ are defined by the Lie algebra structure of $\fg$. Furthermore, the Lie bracket of $\fg$ descends to an action of $\fg/[\fg,\fg]$ on  $[\fg,\fg]/\fz$. If we choose a vector $l\in \fg/[\fg,\fg]$, $l\not=0$, then the eigenvalues of the action of $l$ on $[\fg,\fg]/\fz$ are $\pm t, \pm it\mu_1,\dots,\pm it \mu_q$ for some $t\in\RR$.  Hence, $\mu_1, \dots,\mu_q$ are defined up to signs and up to ordering only by the Lie algebra structure of $\fg$. In the same way we prove the assertion for $\fosc^2_{q}(\mu)$ and $\fosc^2_{q}(\mu')$. Furthermore, the argument shows that $\mu(\RR^2)=\mu'(\RR^2)\ \mod  S_{q}\ltimes(\ZZ_{2})^{q}$ holds if $\fd_{q}(\mu)$ and $\fd_{q}(\mu')$ are isomorphic. Now suppose that $\mu(\RR^2)=\mu'(\RR^2)\ \mod  S_{q}\ltimes(\ZZ_{2})^{q}$. We show that the Lie algebras $\fd_{q}(\mu)$ and $\fd_{q}(\mu')$ are isomorphic. By Prop.~\ref{P1}, we may assume that  $\mu'=s\mu$ for some real number $s\not=0$. The linear map $\phi: \fd_{q}(\mu)\to \fd_{q}(s\mu)$, defined by $\phi(\fz)=\fz, \phi(\fa)=\fa, \phi(\fl)=\fl$ and
$$ \phi|_\fz=(1/s^3)\cdot\Id_\fz,\ \phi|_\fa=(1/s^2)\cdot\Id_\fa,\  \phi|_\fl=(1/s)\cdot\Id_\fl$$
is an isomorphism.
\qed

Propositions \ref{P1} and \ref{P2} imply a classification (up to Lie group isomorphisms) of simply-connected Lie groups that admit a bi-invariant metric of index 2.  In the following subsections we describe these groups. 
\subsection{Lie groups with one-dimensional centre}\label{S2.1}
We want to describe the simply-connected group integrating $\fosc_{1,q}(1,\mu)$. Consider $\fl=\RR$ and let $(\rho,\fa,\ip_\fa)$ be as defined in item (i) in the previous subsection.
Then $\omega:=\langle L(\,\cdot\,),\cdot\rangle_\fa$ is a symplectic form on $\fa$ satisfying $\omega(La,a')+\omega(a,La')=0$ for all $a,a'\in\fa$. The ideal
$$\fn_1(\mu):=\fz\oplus \fa=\RR\oplus \fa$$ of $\fosc_{1,q}(1,\mu)$ has centre $\fz$ and the Lie bracket of $a, a'\in \fa$ equals 
 $$[a, a']=\omega(a,a')\in\fz.$$
In particular, it is isomorphic to the $(n-1)$-dimensional Heisenberg algebra for an arbitrary choice of $\mu\in (\RR\setminus\{0\})^q$.  The Lie algebra $\fosc_{1,q}(1,\mu)$ is equal to the semidirect product  
 $$\fosc_{1,q}(1,\mu)=\fn_1(\mu)\rtimes_\rho \fl,$$
where the representation $\rho$ of $\fl$ on $\fn_1(\mu)=\fz\oplus\fa$ is given by $\rho(t)(z,a)=(0,tL(a))$ for $t\in \fl$. We denote by $N_1(\mu)$ the simply-connected $(n-1)$-dimensional Heisenberg group with Lie algebra $\fn_1(\mu)$. It is a central extension of the abelian Lie group $\fa$ by $\fz=\RR$ and the multiplication is given by 
\begin{equation}
(z,a)\cdot (z', a')=(z+z'+\textstyle{\frac12} \omega(a,a'), a+ a') \label{E*Hn}
\end{equation}
for $z,z' \in\frak z$ and $a,a'\in\fa$. The simply-connected Lie group with Lie algebra $\fosc_{1,q}(1,\mu)$ is equal to the semidirect product 
$$\Osc_{1,q}(1,\mu):= N_1(\mu)\rtimes \RR,$$
where $t\in \RR$ acts on $N_1(\mu)$ by $t\cdot (z,a)=(z, e^{tL}(a))$. If $q=0$, then we write $\Osc_{1,0}$ instead of $\Osc_{1,0}(1,\emptyset)$.

\subsection{Even-dimensional groups with two-dimensional centre}\label{Sev}
We describe the simply-connected Lie group associated with $\fosc^2_{q}(\mu)$. Consider $\fl=\RR^2$ and let $\rho$ and $\fa=\fa^{2q}$ be as defined in item (ii) in Subsection \ref{Sqe}. Then $\fz=\fl^*\cong\RR^2$.
We define $\omega\in \Hom(\bigwedge^2 \fa,\fz)$ by $\omega(a,a')(t):=\la \rho(t)a,a'\ra_{\fa}$ for $t\in\fl$. Then  $$\fn_2 (\mu):=\fz\oplus \fa=\RR^2\oplus \fa$$ is an ideal of $\fosc^2_{q}(\mu)$ with centre $\fz$ and 
 $$[a,a']=\omega(a,a')\in \fz$$ 
for $a,a'\in \fa$ and we obtain
$$\fosc^2_{q}(\mu)=\fn_2(\mu)\rtimes \fl,$$
where $\fl$ acts trivially on $\fz$ and by $\rho$ on $\fa$. 

Let $N_2(\mu)$ denote the simply-connected Lie group with Lie algebra $\fn_2(\mu)$. This group is a central extension of the abelian group $\fa$ by $\fz \cong \RR^2$ and the multiplication is given by 
\begin{equation}
(z,a)\cdot (z', a')=(z+z'+\textstyle{\frac12} \omega(a,a'), a+ a') \label{E*H2}
\end{equation}
for $z,z' \in\frak z$ and $a,a'\in\fa$. The simply-connected Lie group with Lie algebra $\fosc^2_{q}(\mu)$ is equal to the semidirect product 
$$ \Osc^2_q(\mu):= N_2(\mu)\rtimes \RR^2,$$
where $t\in \RR^2$ acts on $N_2(\mu)$ by $t\cdot (z,a)=(z, e^{\rho(t)}(a))$. 
\subsection{Odd-dimensional groups with two-dimensional centre}
Finally, we describe the simply-connected Lie group associated with $\fd_q(\mu)$. Recall that in this case $\fl=\RR^2$, $\fz=\fl^*$ and $\fa=\fa^{2q}\oplus\fa_0$, where $\fa_0=\RR$. Let $\omega\in \Hom(\bigwedge^2 \fa,\fz)$ be defined as in Subsection \ref{Sev}. We consider again $\fn_2 (\mu)=\fz\oplus \fa^{2q}$. The Lie algebra $\fd_q(\mu)$ is an extension of $\fl$ by the direct sum of Lie algebras $\fn_2 (\mu)\oplus\fa_0$. The Lie bracket between two elements of $\fl$ and between an element of $\fl$ and one of $\fn_2 (\mu)\oplus\fa_0$ is given by (\ref{c1}) and (\ref{c2}). In order to keep the notation as simple as possible we put $t^\flat:=\alpha(t,\cdot) \in \fl^*$ for $t\in\fl$. For $s\in\fa_0$ and $t\in\fl$, by $st^\flat\in \fz$ we understand the element $t^\flat$ multiplied by the real number $s$. For $t,t'\in \fl$, $z\in\fz$, $a\in\fa^{2q}$ and $s\in\fa_0$, we then have 
\begin{eqnarray*}
[t,z+a+s]&=&-st^\flat+\rho(t)(a)\ \in\ \fz\oplus\fa^{2q}= \fn_2 (\mu),\\
{}[t,t'] &=& \alpha(t,t') \ \in\ \fa_0.
\end{eqnarray*}
Let $D_q(\mu)$  denote the corresponding simply-connected Lie group. For $q=0$ we write $D_0$ instead of $D_0(\emptyset)$.

For $z\in\fz$, $a\in\fa^{2q}$, $s\in\fa_0$ and $t\in \fl$, we put 
$$l(t):=\exp(t)\in D_q(\mu),\quad h(z,a,s):=\exp(z+a+s)\in D_q(\mu).$$
Since $D_q(\mu)$ is an extension of the simply-connected nilpotent Lie group $N_2(\mu)\times \RR$ by $\RR^2$, the map
$$\fd_q(\mu)\ni z+a+s+t\longmapsto h(z,a,s)\cdot l(t)\in D_q(\mu)$$
is a diffeomorphism. 
\begin{lm} For $t,t'\in\fl$ and $(z,a,s)\in\fn_2(\mu)\oplus \fa_0$, we have
\begin{eqnarray}
h(z,a,s)\cdot h(\hat z,\hat a,\hat s)&=&\textstyle{h(z+\hat z+\frac12\omega(a,\hat a),a+\hat a,s+\hat s)  }\label{Ehh}\\
l(t)\cdot l(\dt) &=& h\left(\textstyle{-\frac13\alpha(t,\dt)(t+\frac12 \dt)^\flat,0,\frac12 \alpha(t,\dt)}\right) \cdot l(t+\dt)   ,      \label{Ell}\\
l(t)\cdot h(z,a,s)\cdot l(t)^{-1}&=&h(z-s t^\flat,e^{\rho(t)}(a),s). \label{Elhl}
\end{eqnarray}
\end{lm}
\proof The exponential map $\exp:\fd_q(\mu)\rightarrow D_q(\mu)$ restricted to $\fn_2(\mu)\oplus \fa_0$ equals $\exp:\fn_2(\mu)\oplus \fa_0\rightarrow N_2(\mu)\times\fa_0$, which gives Eq.\,(\ref{Ehh}).  

Since $l(t)\cdot l(\dt) = \exp(t)\cdot \exp(\dt )$, the Baker-Campbell-Hausdorff formula gives
\begin{eqnarray*}
l(t)\cdot  l(\dt)&=&\exp (t+\dt+\textstyle{ \frac12 [t,\dt]+\frac1{12}[t-\dt,[t,\dt]]}\,)\\
&=& \exp\left(-\textstyle{\frac 1{12} \alpha(t,\dt)(t-\dt)^\flat+\frac12 }\alpha(t,\dt) + t+\dt\,\right).
\end{eqnarray*}
Analogously, 
\begin{eqnarray*}
\lefteqn{ h\left(\textstyle{-\frac13\alpha(t,\dt)(t+\frac12 \dt)^\flat,0,\frac12 \alpha(t,\dt)}\right) \cdot l(t+\dt)   }\\
&&=\exp\left(\textstyle{-\frac13\alpha(t,\dt)(t+\frac12 \dt)^\flat+\frac12 \alpha(t,\dt)}\right) \cdot \exp(t+\dt) \\
&&=\exp\left (\textstyle{-\frac13\alpha(t,\dt)(t+\frac12 \dt)^\flat+\frac12 \alpha(t,\dt)+t+\dt+ \frac12 [\frac12 \alpha(t,\dt) ,t+\dt]}\right) \\
&&=\exp\left (\textstyle{-\frac13\alpha(t,\dt)(t+\frac12 \dt)^\flat+ \frac14\alpha(t,\dt)(t+\dt)^\flat+\frac12 \alpha(t,\dt)+t+\dt }\right) \\
&&= \exp\left(-\textstyle{\frac 1{12} \alpha(t,\dt)(t-\dt)^\flat+\frac12 }\alpha(t,\dt)+ t+\dt\,\right),
\end{eqnarray*}
which proves Eq.\,(\ref{Ell}). 
As for Eq.\,(\ref{Elhl}), 
\begin{eqnarray*}
l(t)\cdot h(z,a,s)\cdot l(t)^{-1}&=&l(t)\cdot \exp(z+a+s)\cdot l(t)^{-1}\\
& =& \exp ( \Ad(l(t))(z+a+s))\ =\ \exp (e^{\ad\, t} (z+a+s))\\
&=&h(z-s t^\flat,e^{\rho(t)}(a),s). 
\end{eqnarray*}
\qed

\section{Existence of lattices}
\subsection{Lie groups with one-dimensional centre}

In this section we discuss the question which simply-connected indecomposable Lie groups with bi-invariant metric of index 2 and with one-dimensional centre have a lattice. We use methods developed in \cite{KOmemo} to give an answer in terms of Salem numbers.

If $p$ is the minimal polynomial over $\QQ$ of an algebraic complex number $z_0\not=1$ with $|z_0| = 1$, then $p$ is self-reciprocal, 
 i.e., it satisfies $p(x)=x^{\deg p}p(1/x)$.
In particular, if $\alpha$ is a root of $p$, then also $1/\alpha$ is a root of $p$. We apply this to the following special kind of polynomials. 
Let $f\in\ZZ[x]$ be a monic polynomial of degree $2k$, $k> 1$, that is  irreducible over $\ZZ$ and has exactly $2k-2$ roots on the unit circle.
Then $f$ is self-reciprocal and the two roots $r_1,r_2$ of modulus different from 1 are real. In particular, $r_1=1/r_2$. Let $F^+_{2k}$ denote the set of all such polynomials~$f$ for which, in addition, the two real roots are positive. This set is related to the notion of Salem numbers. 
We want to quickly recall this notion here, for a general introduction we refer to \cite{Sm}.
By definition, a Salem number is a real algebraic integer $r>0$, all of whose Galois conjugates different from $r$  lie in the closed unit disc $|z| \le 1$, with at least one on its boundary $|z| = 1$. Hence,  $r\in \RR$ is  a Salem number if and only if it is equal to the real root of modulus $>1$ of a polynomial in $F^+_{2k}$ for some $k\ge2$. The number $2k$ is called the degree of $r$. 
Thus, $F^+_{2k}$ consists precisely of the minimal polynomials of Salem numbers of degree $2k$.
It is not hard to prove that 
$$ F_4^+=\{ x^4-ax^3+bx^2-ax+1\mid 2a>|b+2|, \ b\not=2,\,b\not=\pm a+1 \},$$
see also \cite{B1}. Salem numbers of degree 6 are studied, e.g., in \cite{B2}. In the Supplement to \cite{B2} one can find tables listing examples of polynomials contained in $F_6^+$. For instance, $f(x)=x^6-x^4-x^3-x^2+1$ is in $F_6^+$. 

For $k=1$ we put 
$$F_2^+:=\{x^2-ax+1\mid a\in\ZZ,\ a \ge 3\},$$ 
which equals the set of quadratic integer polynomials having real roots $r>1$ and $1/r$.

In the following, the symbol $\{\dots\}_!$ will denote multisets. Elements will be repeated according to their multiplicity.
For $q>0$ and $f\in F_{2\bar q+2}^+$, $q\ge \bar q\ge 0$, we define the set $M_{f, q} \subset \RR^q$ in the following way. Let $r>1,1/r \in\RR$ and $\nu_1^{\pm1},\dots, \nu_{\bar q}^{\pm1}\in S^1$ be the roots of $f$ and put
$$M_{f,q}:=\Big\{(\mu_1,\dots,\mu_q)\in (\RR\setminus\{0\})^q\mid \{r^{ \pm i \mu_1},\dots,r^{ \pm i \mu_q}\}_!= \{\nu_1^{\pm1},\dots,\nu_{\bar q}^{\pm1},\underbrace{1,\dots,1}_{2(q-\bar q)\times}\}_!  \Big\}.$$  
\begin{theo}\label{T1}
The group $\Osc_{1,0}$ admits a lattice. For $q>0$ and $\mu\in(\RR\setminus\{0\})^q$, the group $\Osc_{1,q}(1,\mu)$ admits a lattice if and only if there exists a polynomial $f\in F^+_{2\bar q+2}$ for some $0\le \bar q\le q$  such that $\mu\in M_{f,q}$.
\end{theo}
The proof will rely on the following lemma, which can be found in \cite{KOmemo}, Section 7.
\begin{lm}\label{LKO}
Let $(V,\omega)$ be a symplectic vector space over $\RR$, and let $A\in \Sp(V,\omega)$ be semisimple with characteristic polynomial $f_A$. Then there exists a lattice $\Lambda\subset V$ with $A\Lambda=\Lambda$ and with $\omega(\Lambda\times \Lambda)\subset\ZZ$ if and only if $f_A$ has integer coefficients. 
\end{lm}
{\sl Proof of Theorem~\ref{T1}.}
We will use the following well-known facts:
\begin{itemize} 
\item[-] Let $G$ be a connected solvable Lie group  and $N$ be its maximal connected normal nilpotent Lie subgroup. If $\Gamma$ is a lattice in $G$, then $\Gamma \cap N$ is a lattice in $N$ (see \cite{Ra72}, Cor.\,3.5.).
\item[-] If $\Gamma_0$ is a lattice in the simply-connected nilpotent Lie group $N$, then the image of $\Gamma_0$ in the factor group $N/[N,N]$ is also a lattice (see {\rm \cite{Mal}}, Thm.\,4). 
\end{itemize}
Suppose that $\Gamma$ is a lattice in $\Osc_{1,q}(1,\mu)$, $q>0$. Let $N:=N_1(\mu)\subset \Osc_{1,q}(1,\mu)$ be as defined in Subsection~\ref{S2.1}. Then $N$ is the maximal connected 
normal
nilpotent Lie subgroup of $\Osc_{1,q}(1,\mu)$. Hence $\Gamma_0:=\Gamma \cap N$ is a lattice in $N$.  Furthermore, the image $\Lambda$ of $\Gamma_0$ under the projection $N\to N/[N,N]\cong \fa$ is a lattice in $\fa$. 
Choose an element $(z_0,a_0,t_0)\in \Gamma$ such that $t_0>0$. Then $e^{t_0L}(\Lambda)=\Lambda$ since $\Gamma_0$ is invariant under conjugation by $(z_0,a_0,t_0)$. Choose $n\in \NN$ such that all eigenvalues of $(e^{t_0L})^n=e^{nt_0L}$ that are roots of unity are equal to one. Put $t':=nt_0$ and $A:=e^{t'L}$.  
Then $A\Lambda=\Lambda$,
hence the characteristic polynomial $f_A$ of $A$ has integer coefficients.
 By construction, $e^{\pm t'},e^{\pm it'\mu_1},\dots,e^{\pm it'\mu_q}$ are the eigenvalues of $A$, where  $e^{t'}>1$ and each of the eigenvalues $e^{\pm it'\mu_1},\dots,e^{\pm it'\mu_q}$ either is equal to 1 or is not a root of unity. Consequently, $f_A(x)=f(x)\cdot  (x-1)^{2(q-\bar q)}$ for some monic polynomial  $f\in\ZZ[x]$ of degree $2\bar q+2$ whose roots are not roots of unity. The roots of $f$, except  $e^{\pm t'}$, all have modulus 1. Hence $f$ is irreducible 
over $\ZZ$.
 Indeed, assume that $f=f_1f_2$ for polynomials $f_1,f_2\in\ZZ[x]$ and that  $e^{t'}$ is a root of $f_1$. Then all roots of the monic polynomial $f_2\in\ZZ[x]$  have modulus at most 1 (and are non-zero). 
Hence all these roots are roots of unity by a well-known result of Kronecker 
\cite{Kr}.
 But this contradicts our assumption on the roots of $f$.  Hence $f$ is irreducible and we conclude $f\in F^+_{2\bar q+2}$. The real eigenvalues of $f$ are $r:=e^{t'}>1$ and $1/r$. Hence $\mu\in M_{f, q}$. 

Now suppose that $q>0$ and take $\mu\in M_{f, q}$. 
Let $r>1$ be one of the real roots of $f$. Then $r=e^{t'}$ for some $t'>0$. By assumption the characteristic polynomial of $A:=e^{t'L}$ equals $f_A(x)=f(x)\cdot (x-1)^{2(q-\bar q)}$. Since also $A\in\Sp(\fa,\omega)$, Lemma~\ref{LKO} ensures the existence of a lattice $\Lambda$ of $\fa$ satisfying $A\Lambda=\Lambda$ and $\omega(\Lambda\times\Lambda)\subset \ZZ$. Then $\frac12\cdot\ZZ\times\Lambda\times t'\cdot\ZZ\subset \fz\times \fa\times\fl$ is a lattice in $\Osc_{1,q}(1,\mu)$. 

If $q=0$, choose  $t'>0$ such that $e^{\pm t'}$ are the roots of a polynomial in $F^+_2$. Then put $A:=e^{t'L}$ and proceed as above. Hence $\Osc_{1,0}$ has a lattice.
\qed
\begin{re}\label{salem}{\rm
By definition, $(\mu_1,\dots,\mu_q)\in (\RR\setminus\{0\})^q$ is in $M_{f,q}$ if and only if
\begin{eqnarray*}
\mu_j&\in& \pm \frac{s_j}{\ln r}+\frac{2\pi}{\ln r}\cdot \ZZ,\quad j=1,\dots,\bar q,\\
\mu_j&\in & \frac{2\pi}{\ln r} \cdot (\ZZ\setminus\{0\}),\quad j=\bar q+1,\dots,q,
\end{eqnarray*}
up to permutation of $\mu_1,\dots,\mu_q$, where $r^{\pm1}$ and $\nu_j^{\pm1}=e^{\pm is_j}$, $j=1,\dots,\bar q$, are the zeroes of $f$.
In particular, $M_{f,q}$ is countable.

For $f\in F_2^+$ with roots $r^{\pm1}$, we have $$M_{f,q}=\frac {2\pi}{\ln r} \cdot (\ZZ\setminus \{0\})^q.$$
The root $r$ of $f(x)=x^2-ax+1\in F^+_2$ with $r>1$ satisfies $\displaystyle{\lim_{a\to\infty} r=\infty}$, thus 
$\bigcup_{f\in F_2^+} M_{f,q}$ is dense in $\RR^q$.
}\end{re}

\begin{co}
The set of those $\mu\in(\RR\setminus\{0\})^q$ for which $\Osc_{1,q}(1,\mu)$ has a lattice is countable and dense in $\RR^q$.
\end{co}

Before we give an example, we want to review some basic properties of Salem numbers, see e.g., \cite{Sa,Sm} for further information. If $r$ is a Salem number of degree $d$, then so is $r^k$ for all $k\in \NN_{>0}$. Suppose that $r$ and $r'$ are Salem numbers and that  $r'$ is an element of the field extension $\QQ(r)$. Then $\QQ(r') = \QQ(r)$.  If $r$ is a Salem number, then there is a Salem number $\tau \in \QQ(r)$ such that the set of Salem numbers in $\QQ(r)$ consists of the powers of $\tau$. We define an equivalence relation on the set of Salem numbers:
\begin{equation}\label{equiv} 
r_1 \sim r_2\ :\Leftrightarrow \exists \ k_1, k_2\in \NN_{>0}:r_1^{k_1}=r_2^{k_2}.
\end{equation}
By the above mentioned properties of Salem numbers, $r_1\sim r_2$ holds if and only if $\QQ(r_1)=\QQ(r_2)$, which is equivalent to the existence of a Salem number $r$ such that $r_1=r^{l_1}$ and $r_2=r^{l_2}$ for some $l_1, l_2\in\NN_{>0}$.
\begin{ex} {\rm We consider the case  $q=1$, i.e., generalised oscillator groups $\Osc_{1,1}(1,\mu)$, $\mu>0$. Our aim is to determine the set ${\cal M}$ of parameters $\mu>0$ for which $\Osc_{1,1}(1,\mu)$ admits a lattice. By Theorem~\ref{T1}, we know that ${\cal M}={\cal M}_2\cup{\cal M}_4$, where
$${\cal M}_k=\bigcup _{f\in F_k^+} (M_{f,1}\cap \RR_{>0}),\quad k=2,4.$$ 
We want to describe the sets ${\cal M}_2$ and ${\cal M}_4$ explicitly. We denote by $T_2$ the set of real numbers $r>1$ with $f(r)=0$ for some $f\in F_2^+$ and by $T_4$ the set of Salem numbers of degree $4$. By Remark~\ref{salem} we have
\begin{equation}\label{EF2}
{\cal M}_2=\bigcup_{r\in T_2}\frac{2\pi}{\ln r} \NN_{>0}.
\end{equation}
This is not a disjoint union. We want to rewrite it in a way that avoids multiple listing of elements. We put $${\cal M}_{2,r}:= \frac{\pi}{\ln r} \QQ_{>0}$$ for arbitrary $r>0$. For every $r\in T_2$, also each power $r^k$, $k\in \NN_{>0}$, is in $T_2$. Hence,  Eq.\,(\ref{EF2}) implies
${\cal M}_2=\bigcup_{r\in T_2} {\cal M}_{2,r}$. Thus we have a disjoint union 
\begin{equation}\label{EF2a}
{\cal M}_2=\bigcup_{r\in\, T_2/_\sim} {\cal M}_{2,r}, 
\end{equation}
where $r_1\sim r_2$ is also defined as in (\ref{equiv}). The map $T_2\ni r \mapsto \QQ(r)\subset\RR$ descends to a bijection from $T_2/_\sim$ to the set of real quadratic number fields. Furthermore, the equivalence class of $r$ is contained in $K:=\QQ(r)$ and equals the set of units $u>1$ in ${\cal O}_K$ that have norm  $N(u)=1$. All units are powers of the fundamental one. Hence, ${\cal M}_{2,r}={\cal M}_{2,u}$ for every unit $u>1$ in ${\cal O}_K$ and we can write (\ref{EF2a}) also as
$${\cal M}_2=\bigcup_D  {\cal M}_{2,u_D},$$
where the (disjoint) union runs over all square free integers $D\ge 2$ and $u_D>1$ is any unit in the ring of integers in $\QQ(\sqrt D)$.

Next we compute ${\cal M}_4$. Take $f\in F_4^+$ and let $r>1$, $1/r$ and $e^{\pm is}$ be the roots of $f$. Then
\begin{equation}\label{e1}
M_{f,1}\cap \RR_{>0}= \{\mu \in\RR \mid \exists k\in\ZZ: \ln(r)\cdot \mu =|s+2k\pi| \}.
\end{equation}
For $r\in T_4$, we define
$${\cal M}_{4,r}:= \{ |s+q\pi|/\ln(r) \mid q\in\QQ\},$$
where $e^{\pm is}$ are the complex Galois conjugates of $r$.
For $r\in T_4$ and $l\in\NN_{>0}$, also $r^l$ is in $T_4$ and  $e^{\pm ils}$ are Galois conjugates of $r^l$ if $e^{\pm is}$ are Galois conjugates of $r$, hence (\ref{e1}) gives ${\cal M}_4=\bigcup_{r\in T_4} {\cal M}_{4,r}$. Since $r_1\sim r_2$ implies ${\cal M}_{4,r_1}={\cal M}_{4,r_2}$, we obtain 
\begin{equation}\label{EF2b}
{\cal M}_4=\bigcup_{r\in\, T_4/_\sim} {\cal M}_{4,r}.
\end{equation}
The question whether this union is disjoint is
related to the {\sl four exponentials conjecture}, which can be formulated as follows:
Let $(\lambda_{ij})$ be a $(2\times 2)$-matrix of complex numbers such that $e^{\lambda_{ij}}$ is algebraic for all $i,j$ and
such that rows and columns are linearly independent over $\QQ$. Then its rows (and hence the columns) are linearly independent over $\CC$  \cite{Wa}.

Assume that $\mu=|s_1+q_1\pi|/\ln(r_1)=|s_2+q_2\pi|/\ln(r_2)$ for non-equivalent Salem numbers $r_1,r_2\in T_4$ and $q_1,q_2\in\QQ$. Put $\lambda_{11}:=\ln(r_1)$, $\lambda_{12}:=i|s_1+q_1\pi|$, $\lambda_{21}:=\ln(r_2)$ and $\lambda_{22}:=i|s_2+q_2\pi|$. Then all numbers $e^{\lambda_{ij}}$ are algebraic since they are (up to multiplication by  roots of unity) roots of a polynomial in $F_4^+$. Moreover, $\lambda_{11}$ and $\lambda_{12}$ are obviously independent over $\QQ$ and  $\lambda_{11}$ and $\lambda_{21}$ are independent over $\QQ$ since $r_1$ and $r_2$ are not equivalent. Hence the assumptions in the four exponentials conjecture are satisfied. Consequently, if the conjecture is true, then we get a contradiction to our assumption, which implies that the union in (\ref{EF2b}) is disjoint.

Furthermore, if the four exponentials conjecture is true, then also ${\cal M}_2$ and ${\cal M}_4$ are disjoint.  Indeed, assume that $\mu=q_1\pi/\ln(r_1)=|s_2+q_2\pi|/\ln(r_2)$ for $r_1\in T_2$, $r_2\in T_4$ and $q_1,q_2\in\QQ$. Then $\lambda_{11}:=i q_1\pi$, $\lambda_{12}:=\ln (r_1)$, $\lambda_{21}:=i|s_2+\pi q_2|$ and $\lambda_{22}:=\ln (r_2)$ satisfy the assumptions of the conjecture since $e^{\pm is_2}$ are not  roots of unity. Again the conjecture would give a contradiction to our assumption.

To finish the computation of ${\cal M}$ it remains do determine the sets $T_2$ and $T_4$. Obviously, $T_2=\{\frac 12(a+\sqrt{a^2-4})\mid a\in\ZZ,\, a\ge 3\}.$
Also Salem numbers of degree four and their Galois conjugates can be computed by explicit formulas. Indeed, consider the polynomial
$f(x)=x^4-ax^3+bx^2-ax+1\in F_4^+$ with zeroes $r>1$, $1/r$ and $\nu^{\pm1}=e^{\pm is}$. Then 
$$
f(x)=(x-r)(x-1/r)(x-\nu)(x-\bar \nu)=(x^2-t_1x+1)(x^2-t_2x+1),
$$
where $t_1=r+1/r$ and $t_2=2\Re(\nu)$. Hence $t_1+t_2=a$ and $t_1t_2=b-2$, thus $t_1$ and $t_2$ are the roots of the quadratic polynomial $t^2-at+b-2$. Since $t_1>2>t_2$, we get
$$\textstyle{ t_1=\frac a2 + \sqrt{\frac {a^2} 4 -b+2}, \quad t_2=\frac a2 - \sqrt{\frac {a^2} 4 -b+2}}$$
and, finally,
$$\textstyle{r=\frac12\Big( t_1+\sqrt{t_1^2-4}\Big),\quad s=\arccos (t_2/2)}.$$
}\end{ex}
In the remaining part of this subsection, we want to study the classification of abstract commensurability classes of lattices in indecomposable Lie groups that admit a bi-invariant metric of index 2 and have a one-dimensional centre. Recall that two groups $G_1$ and $G_2$ are abstractly commensurable if there are subgroups $H_1\subset G_1$ and $H_2\subset G_2$ of finite index such that $H_1$ is isomorphic to $H_2$. We will see that the equivalence classes are represented by a special kind of extensions of the discrete Heisenberg group by $\ZZ$.
Since there are different notions of discrete Heisenberg groups, we want to explain which one we will use here. Let $e_1,\dots,e_{2m}\in \ZZ^{2m}\subset \RR^{2m}$ denote the standard basis of $\RR^{2m}$. The standard symplectic form $\omega_0$ on $\RR^{2m}$ is defined by $\omega_0(e_{2i-1},e_{2i})=1$ and $\{e_{2i-1},e_{2i}\}\perp \{e_{2j-1},e_{2j}\}$ for $i\not=j$. In the following, by a symplectic map we always mean a map that is symplectic with respect to $\omega_0$. We define a group structure on $\frac12 \ZZ\times \ZZ^{2m}$ by $(z,v)\cdot (z', v')=(z+z'+\textstyle{\frac12} \omega_0(v,v'), v+ v')$.
By a discrete Heisenberg group ${\rm H}_m(\ZZ)$ we understand the subgroup of $\frac12 \ZZ\times \ZZ^{2m}$ generated by elements $\alpha_i:=(0,e_{2i-1}), \beta_{i}:=(0,e_{2i}),\gamma:=(1,0)$. The generators satisfy the relations $[\alpha_i,\beta_i]=\gamma$, $i=1,\dots,m$ and the commutators for all other pairs of generators are trivial. In the following we denote the centre of ${\rm H}_m(\ZZ)$ by $Z$. 

If $a:\ZZ\to \Aut({\rm H}_m(\ZZ))$ is a homomorphism, we can form the semi-direct product ${\rm H}_m(\ZZ)\rtimes _a \ZZ$. The automorphism $\bar a(1)$ induced by $a(1)$ on ${\rm H}_m(\ZZ)/Z\cong\ZZ^{2m}$ is symplectic or anti-symplectic. If we have two homomorphisms $a,a':\ZZ\to \Aut({\rm H}_m(\ZZ))$ for which the induced automorphisms $\bar a(1)$ and $\bar a'(1)$ on $\ZZ^{2m}$  coincide, then the semi-direct products with respect to $a$ and $a'$ are isomorphic. Furthermore, for a given $A\in \Sp(2m,\ZZ)$, we can find a homomorphism $a:\ZZ\to \Aut({\rm H}_m(\ZZ))$ such that the map $\bar a(1)$ induced by $a(1)$ on $\ZZ^{2m}$ equals $A$. Hence, we may define $\Gamma(A):={\rm H}_m(\ZZ)\rtimes _a \ZZ$, where $a:\ZZ\to \Aut({\rm H}_m(\ZZ))$ is any homomorphism such that $\bar a(1)=A$. 
 
Let $\cA_q$ denote the set of abstract commensurability classes of discrete groups that are lattices in an oscillator group $\Osc_{1,q}(1,\mu)$ for some fixed $q>0$.

\begin{pr} 
The elements of $\cA_q$ are (non-uniquely) represented by the groups $\Gamma(A)$, where $A=\diag(I_{2(q-\bar q)}, A^\circ)\in \Sp(2q+2,\ZZ)$, $0\le \bar q\le q$, is a block diagonal matrix and the characteristic polynomial of $A^\circ$ is in $F^+_{2\bar q+2}$. 

For $A_i=\diag(I_{2(q-\bar q_i)}, A_i^\circ)$, $i=1,2$, the groups $\Gamma(A_1)$ and $\Gamma(A_2)$ represent the same element of $\cA_q$ if and only if $\bar q_1=\bar q_2=:\bar q$ and if there exist 
numbers $n_1,n_2\in\ZZ_{\not=0}$ such that $(A_1^\circ)^{n_1}$ and $(A_2^\circ)^{n_2}$ are conjugate in $\GSp(2\bar q+2,\QQ)$.
\end{pr}
\proof
Let $\Gamma$ be a lattice in $\Osc_{1,q}(1,\mu)$, $q>0$. As above, we denote by $\Gamma_0$ its intersection with $N=N_1(\mu)$ and by $\Lambda$ the projection of $\Gamma_0$ to $N/[N,N]$. In the proof of Theorem~\ref{T1} we have seen that there exists an element $\delta=(z',a',t')\in\Gamma$ with $t'>0$ such that all eigenvalues of $A:=e^{t'L}$ that are roots of unity are equal to one. The linear map $A$ on $\fa=\RR^{1,1}\oplus \RR^{2q}$ is semisimple and has characteristic polynomial $f(x)\cdot (x-1)^{2(q-\bar q)}$, where $f\in F^+_{2\bar q+2}$. Moreover, it is symplectic with respect to $\omega= \langle L(\,\cdot\,),\cdot\rangle_\fa$. The lattice generated by $\Gamma_0$ and $\delta$ has finite index in $\Gamma$. Hence we may assume that it is equal to $\Gamma$. We choose a basis of $\RR^{2q+2}$ consisting of generators of $\Lambda$. Then the matrix of $A$ with respect to this basis is integral. We consider $A$ as a linear map on $\QQ^{2q+2}:=\Span_{\Bbb Q}\Lambda$. Then $\QQ^{2q+2}$ is spanned by the orthogonal sum of the kernel and the image of $A-I_{2q+2}$. This sum is orthogonal with respect to $\omega$. We can choose an integral basis $v_1,\dots,v_{2(q-\bar q)}\in\ZZ^{2q+2}$ of $\Ker(A-I_{2q+2})$ and an (integral) basis $w_1,\dots,w_{2\bar q+2}$ of $\Im(A-I_{2q+2})$. Since $\Gamma_0$ is a lattice in $N$, the lattices in $\Ker(A-I_{2q+2})$ and $\Im(A-I_{2q+2})$ spanned by $v_1,\dots,v_{2(q-\bar q)}$ and $w_1,\dots,w_{2\bar q+2}$, respectively, are Heisenberg lattices (with respect to $\omega$) in the sense of \cite{T}. By \cite{T}, Lemma 1.8, we find generators $\bar \alpha_1,\bar\beta_1,\dots,\bar\alpha_{q-\bar q},\bar\beta_{q-\bar q}$ of $\Span_{\Bbb Z}\{v_1,\dots,v_{2(q-\bar q)}\}$ and $\bar \alpha_{q-\bar q+1},\bar\beta_{q-\bar q+1},\dots,\bar\alpha_{q+1},\bar\beta_{q+1}$ of $\Span_{\Bbb Z}\{w_1,\dots,w_{2\bar q+2}\}$ such that $\{\bar\alpha_i,\bar\beta_i\}\perp\{\bar \alpha_j,\bar\beta_j\}$ with respect to $\omega$ for $i\not=j$. Moreover, if $c$ is a generator of the centre $Z(\Lambda)\subset\fz=\RR$ of $\Lambda$, then $\omega(\bar\alpha_i,\bar \beta_i)=l_ic$ for integers $l_i$, $i=1,\dots,q+1$. Multiplying the generators by appropriate integers, we may assume that $\omega(\bar\alpha_i,\bar \beta_i)=lc$ for some $l\in\NN_{>0}$ and all $i=1,\dots,q+1$. Now choose elements $\alpha_i,\beta_i$, $i=1,\dots,q+1$, of $\Gamma_0$ whose projections to $\Lambda$ are equal to $\bar\alpha_i$ and $\bar\beta_i$, respectively, and put $\gamma:=(lc,0)\in N$. Then $\alpha_i,\beta_i$, $i=1,\dots, q+1$ generate a finite index subgroup $H$ of $\Gamma_0$, which is isomorphic to $\dH$. Let $\delta\in\Gamma$ be defined as above such that $\Gamma$ is generated bei $\Gamma_0$ and $\delta$. Then $\delta$ acts by conjugation on $\Gamma_0$. Since $H\subset\Gamma_0$ has finite index and since for any finitely generated group the set of its subgroups of fixed finite index is finite (see e.g. \cite{H}), we see that $\delta^k(H)=H$ for a suitable $k\in\NN_{>0}$. Hence the finite index subgroup of $\Gamma$ generated by $H$ and $\delta^k$ is isomorphic to $\Gamma(\tilde A)$, where $\tilde A=\diag(I_{2(q-\bar q)}, \tilde A^\circ)\in \Sp(q+1,\ZZ)$ and the characteristic polynomial $\tilde f$ of $\tilde A^\circ$ is in  $F^+_{2\bar q+2}$. More exactly, if $f$ is the minimal polynomial of $r>0$, then $\tilde f$ is the minimal polynomial of $r^k$.
 
Before proving the second statement, we verify the following property of homomorphisms between groups of the form $\Gamma(A)$. Let $A=\diag(I_{2(q-\bar q')},A^\circ)\in\Sp(2q+2,\ZZ)$ and $B=\diag(I_{2(q-\bar q)},B^\circ)\in\Sp(2q+2,\ZZ)$ be such that the characteristic polynomials of $A^\circ$ and $B^\circ$ are in $F^+_{2+2\bar q'}$ and $F^+_{2+2\bar q}$, respectively.  Then each injective homomorphism $F:\Gamma(B)\rightarrow \Gamma(A)$ maps $\dH\subset\Gamma(B)$ to $\dH\subset\Gamma(A)$. Indeed, the elements $\alpha_{q-\bar q+1},\dots,\alpha_{q+1}$,  $\beta_{q-\bar q+1},\dots,\beta_{q+1}$ and $\gamma$ generate the commutator group $[\Gamma(B),\Gamma(B)]$. Hence the images of these elements are in $[\Gamma(A),\Gamma(A)]\subset \dH$. Now take $\xi\in\{\alpha_j,\beta_j\mid j=1,\dots, q-\bar q\}\subset \Gamma(B)$. Then $\xi$ commutes with $\alpha_{q+1}$, thus $F(\xi)$ commutes with $F(\alpha_{q+1})$. Since $F(\alpha_{q+1})$ belongs to $[\Gamma(A),\Gamma(A)]$, the projection $\overline{F(\alpha_{q+1})}$ of $F(\alpha_{q+1})$ to $\ZZ^{2q+2}$ is in the image of 
$A-I_{2q+2}$. Assume that $F(\xi)$ is not in $\dH$ but contains a power $\delta^k$ of the generator $\delta:=1$ of the $\ZZ$-factor  in $\Gamma(A)=\dH\rtimes \ZZ$. Then the image of $A-I_{2q+2}$ would contain a non-trivial element that is fixed under $A^k$, namely $\overline{F(\alpha_{q+1})}$. But this is impossible since no eigenvalue of the restriction of $A$ to the image of $A-I_{q+1}$ is a root of unity. Thus also $F(\xi)$ is in $\dH$.

Now suppose that $\Gamma(A_1)$ and $\Gamma(A_2)$ are abstractly commensurable, where $$A_i=\diag(I_{2(q-\bar q_i)},A_i^\circ)\in\Sp(2q+2,\ZZ),\ i=1,2,$$ and the characteristic polynomial of $A_i^\circ$ belongs to $F^+_{2\bar q_i+2}$. Let $\Gamma$ be a group that can be embedded as a finite index subgroup into both $\Gamma(A_1)$ and $\Gamma(A_2)$. Then the image of $\Gamma$ is also a lattice in a generalised oscillator group $\Osc_{1,q}(1,\mu)$. Hence $\Gamma$ contains a finite index subgroup of the form $\Gamma(B)$ for some $B=\diag(I_{2(q-\bar q')},B^\circ)\in\Sp(2q+2,\ZZ)$, where the characteristic polynomial of $B^\circ$ is in $F^+_{2+2\bar q'}$. Hence, there exist injective homomorphisms $F_i:\Gamma(B)\rightarrow\Gamma(A_i)$, $i=1,2$. As we have seen above, $F_i(\dH)\subset \dH$, $i=1,2$. Let $S_i$ denote the map induced by $F_i |_{{\rm H}_{q+1}(\Bbb Z)}$ on $\ZZ^{2q+2}$. Then $S_i^*\omega_0=m_i\omega_0$ and $S_iB=A_i^{n_i}S_i$ for appropriate $m_i,n_i\in\ZZ_{\not=0}$ since $F_1$ and $F_2$ are homomorphisms. Considered as linear maps on $\QQ^{2q+2}$, $S_1$ and $S_2$ are invertible. If we put $S:=S_2S_1^{-1}$, then $S^*\omega_0=\frac{m_2}{m_1}\omega_0$ and $A_2^{n_2}=SA_1^{n_1}S^{-1}$. In particular, the multiplicities of the eigenvalue 1 of $A_1^{n_1}$ and $A_2^{n_2}$ are equal, thus those of $A_1$ and $A_2$ are equal (where we again use that 1 is the only eigenvalue which is a root of unity). This implies $\bar q_1=\bar q_2=:\bar q$.  Furthermore, it follows, that $S$ has block diagonal form $S=\diag(S',S'')$, where $S'\in \GSp(2(q-\bar q),\QQ)$ and $S'' \in \GSp(2\bar q+2,\QQ)$. In particular, $(A_2^\circ)^{n_2}=S''(A_1^\circ)^{n_1}(S'')^{-1}$. 
 
Conversely, suppose that there is a map $S'' \in \GSp(2\bar q+2,\QQ)$ such that $S''(A_1^\circ)^{n_1}=(A_2^\circ)^{n_2}S''$. Multiplying by an appropriate integer, we may assume that $S''$ is integral and that $(S'')^*\omega_0=m\omega_0$ for some $m\in\ZZ_{\not=0}$, where we denoted the restriction of $\omega_0$ to $\QQ^{2\bar q+2}$ also by $\omega_0$. We choose an integral matrix $S'\in \GL(2(q-\bar q),\QQ)$ such that also $(S')^*\omega_0=m\omega_0$ on $\QQ^{2(q-\bar q)}$, where we again denote the restriction of $\omega_0$ also by $\omega$. Then $S:=\diag(S',S'')$ is an integral matrix satisfying $S^*\omega_0=m\omega_0$ and $SA_1^{n_1}=A_2^{n_2}S$. We extend $S$ to an injective homomorphism $f:\dH\rightarrow \dH$. Let $a_2:\ZZ\to \Aut(\dH)$ be such that the map $\bar a_2(1)$ induced by $a_2(1)$ on $\ZZ^{2m}$ equals $A_2$, i.e., $\Gamma(A_2)=\dH\rtimes_{a_2}\ZZ$. We choose $k\in\NN_{>0}$ such that the image of $f$ is invariant under $a_2(kn_2)\in  \Aut(\dH)$, where we again use the above cited result from \cite{H}. We put $a:=(a_2)^{kn_2}$ and $b:=f^{-1}af$. We define an extension of $f$ to an isomorphism 
$$F: \dH \rtimes_b \ZZ \longrightarrow f(\dH)\rtimes_a\ZZ\cong  f(\dH)\rtimes_{a_2} kn_2\ZZ \subset \Gamma(A_2)$$ 
which is the identity on the $\ZZ$-factor. Let $B$ be the map induced by $b(1)$ on $\ZZ^{2q+2}$. Then $B=S^{-1}A_2^{kn_2}S=A_1^{kn_1}$. Thus  $\Gamma(B)\cong\Gamma(A_1^{kn_1})$. We obtain $$\Gamma(A_1)\sim \Gamma(A_1^{kn_1})\cong \Gamma(B) \sim\Gamma(A_2),$$
where $\sim$ denotes equivalence under abstract commensurability. \qed
\begin{re}
{\rm 
All (minimal polynomials of) Salem numbers appear in this classification of abstract commensurability casses. More exactly, for every $f\in F^+_{2\bar q+2}$ there exists a map $A^\circ \in \Sp(2\bar q+2,\ZZ)$ such that  the characteristic polynomial of $A^\circ$ equals $f$. In fact, for every monic self-reciprocal polynomial $p$ of degree $2n$ over $\ZZ$ one can find a matrix $M\in \Sp(2n, \ZZ)$ whose characteristic polynomial equals $p$, cf. \cite{Ki}, see also \cite{Ri}.
 
Let us fix a Salem number of degree $2k$ with minimal polynomial $f$ and consider the set $\cal S$ of all elements of $\Sp(2k,\QQ)$ that have characteristic polynomial $f$. Conjugacy classes in $\cal S$ with respect to $\Sp(2k,\QQ)$ can be described, see \cite{Wall}. This description can easily be extended to concugacy with respect to $\GSp(2k,\QQ)$, which leads to an additional $\QQ^\times$-factor. More exactly, let $\QQ(r)\to \QQ(r),\ w\mapsto w^*$ be the automorphism of $\QQ(r):\QQ$ that maps $r$ to $r^{-1}$. The fixed field of this automorphism equals $\QQ(\alpha)$, where $\alpha:=r+r^{-1}$. Then $\cal S$ is in bijection with the quotient of the group of units $\QQ(\alpha)^\times$ by $\{qww^*\mid q\in\QQ^\times, w\in \QQ(r)^\times\}$.  
}\end{re}
\subsection{Groups with two-dimensional centre}
Recall that we here consider only Lie groups with bi-invariant metric which are indecomposable. By Prop.~\ref{P1}, $\Osc^2_q(\mu)$ is indecomposable if and only if $\mu_j\not=0$ for $j=1,\dots,q$ and if the set $\{\mu_1,\dots,\mu_q\}\subset (\RR^2)^*$ is not contained in the 
union of two 1-dimensional subspaces (in particular, only if $q\ge3$). Furthermore, $D_q(\mu)$ is indecomposable if and only if $\mu_j\not=0$ for all $j=1,\dots,q$.
\begin{theo} \label{T2} Let $G$ be one of the groups $\Osc^{2}_q(\mu)$ or $D_q(\mu)$ for some $\mu\in((\RR^2)^*)^q$, $q>0$, and suppose that $G$ is indecomposable. Then $G$ admits a lattice if and only if $\mu_1,\dots,\mu_q$ are contained in a lattice of $(\RR^2)^*$. The group $D_0$ admits a lattice.
\end{theo}
\proof Assume first that $q>0$ and that $\mu_1,\dots,\mu_q$ are contained in a lattice $\Gamma^{\fz}$ of $\fz=(\RR^2)^*$. Let $\sigma_1$ and $\sigma_2$ generate this lattice. Using the standard basis $e_1,\dots,e_q$ of $\CC^{2q}$ we obtain an orthonormal basis $e_1, ie_1, \dots, e_q, ie_q$ of $\fa^{2q}$. The basis vectors span the lattice $\Lambda=\ZZ^{2q}$ in $\fa^{2q}=\RR^{2q}$.  Since
$$
\omega(e_k, ie_k)=\langle \rho(\cdot) e_k,ie_k \rangle_\fa =  \langle \mu_k(\cdot) i e_k,ie_k\rangle_\fa =\mu_k\in \Gamma^{\fz} 
$$
and $\omega(e_k,e_l)=\omega(e_k,ie_l)=\omega(ie_k,ie_l)=0$ for $k\not=l$, we have
\begin{equation}\label{omega}
\omega(\Lambda,\Lambda)\subset  \Gamma^\fz.
\end{equation}
Equation (\ref{E*H2}) now implies that the set $\Gamma^N:=\textstyle{\frac12} \Gamma^\fz\times \Lambda \subset N_2(\mu)$ is a subgroup and hence a lattice in $N_2(\mu)$. Denote by $T_1,T_2$ the basis of $\RR^2$ satisfying $\sigma_i(T_j)=2\pi \delta_{ij}$. Then $T_1,T_2$ span a lattice $\Gamma'$ of $\RR^2$.  The action of $\Gamma'\subset\RR^2$ on $N_2(\mu)$ by conjugation is trivial. Indeed, for $t\in \Gamma'$, we have $\mu_j(t)\in \ZZ \sigma_1(t)+\ZZ \sigma_2(t)\in 2\pi\ZZ$, hence $t$ acts on $N_2(\mu)$ by
$$(z,a)\longmapsto (z,e^{\rho(t)}a)=(z,a) .$$
This shows that $\Gamma^N\rtimes \Gamma'$ is a lattice in $\Osc^{2}_q(\mu)$. 

For $D_q(\mu)$ we have to modify this construction. Let 
$\sigma_1$, $\sigma_2$, $T_1$, $T_2$ and $\Gamma^\fz$, $\Gamma'$,  $\Lambda$ as defined before.
We put $\lambda:=\alpha(T_1,T_2)/\sqrt{2\pi}$ and consider the set
$$\textstyle{\Gamma^\fh:=\frac{1}{6}\lambda^2 \Gamma^{\fz} \times  \lambda \Lambda \times \sqrt{ \frac\pi 2}\lambda \ZZ},$$  
which we consider as a subset of the Lie algebra $\fn_2(\mu)\oplus \fa_0\subset \fd_q(\mu)$.
Let us show that $h(\Gamma^\fh)\cdot l(\Gamma')$
is a subgroup of $D_q(\mu)$. First note that $h(\Gamma^\fh)$ is a subgroup of $D_q(\mu)$. Indeed, by (\ref{omega}),  we have
$$ 
\textstyle{\frac12\omega(\lambda \Lambda,\lambda \Lambda)\subset \frac {\lambda^2}{2}\Gamma^\fz\subset \frac {\lambda^2}{6}\Gamma^\fz},
$$
hence $h(\Gamma^\fh)\cdot h(\Gamma^\fh)\subset h(\Gamma^\fh)$ by (\ref{Ehh}). For $t=n_1T_1+n_2T_2\in \Gamma'$, we have 
$$\mu_j(t)\in (\ZZ \sigma_1+\ZZ \sigma_2)(t)=(\ZZ \sigma_1+\ZZ \sigma_2)(n_1T_1+n_2T_2)\subset 2\pi\ZZ,$$ 
thus $e^{\rho(t)}=\Id_\fa$. Furthermore, for $s\in \sqrt{ \frac\pi 2}\lambda \ZZ $ and for $t\in \ZZ T_1+\ZZ T_2$, we have
\begin{equation}\label{hilf}
t^\flat= \alpha(t,\cdot)=\textstyle{\frac1{2\pi}\alpha(t,T_1)\sigma_1+ \frac1{2\pi}\alpha(t,T_2)\sigma_2}\in \frac\lambda{\sqrt{2\pi}}\Gamma^\fz,
\end{equation} 
thus
$$\textstyle{s t^\flat\in \sqrt{ \frac\pi 2}\lambda  \cdot  \frac\lambda{\sqrt{2\pi}}\Gamma^\fz   = \frac{\lambda^2}{2} \Gamma^\fz\subset \frac{\lambda^2}{6} \Gamma^\fz.}$$
Now, (\ref{Elhl}) implies 
$\textstyle{l(T)\cdot h(\Gamma^\fh)\cdot l(T)^{-1}\subset h(\Gamma^\fh)}$. Finally,
$$\textstyle{\frac13\alpha(t,\dt)(t+\frac12 \dt)^\flat\in \frac13\lambda \sqrt{2\pi} \cdot \frac\lambda{2\sqrt{2\pi}}\Gamma^\fz=\frac{\lambda^2}{6} \Gamma^\fz}$$
for $t,t'\in \Gamma'$. Using (\ref{Ell}), we obtain $l(\Gamma')\cdot l(\Gamma')\in h(\Gamma^\fh)\cdot l(\Gamma')$. Obviously, $\Gamma$ is discrete and cocompact, hence a lattice.

As for $q=0$, recall that $D_0$ is an extension of $\fl=\RR^2$ by the abelian group $\fz\times \fa_0$. In this case, the set 
$$\textstyle{\frac16\ZZ^2}\times \frac 12 \ZZ \times \ZZ^2 \subset \fz\times\fa_0\times \fl$$
is a subgroup of $D_0$. Obviously it is discrete and cocompact.

Now suppose that $\Gamma$ is a lattice in $G$. We can argue as in the proof of Thm.~\ref{T1}. Let $N$ be the maximal connected nilpotent subgroup of $G$, i.e. $N=N_2(\mu)$ if $G=\Osc^2_q(\mu) $ and $N=N_2(\mu)\times\fa_0$ if $G=D_q(\mu)$. Then $\Gamma_0:=\Gamma\cap N$ is a lattice in $N$ and the image $\Lambda$ of $\Gamma_0$ in the factor group $N/[N,N]\cong \fa$ is a lattice in $\fa$.  Choose a basis $T_1,T_2$ of $\fl=\RR^2$ contained in the projection of $\Gamma$ to $\fl$, that is, $T_i=\proj_\fl(\gamma_i$) for $\gamma_i\in\Gamma$, $i=1,2$. Since $\Gamma_0$ is invariant under conjugation by $\gamma_1$ and $\gamma_2$, the lattice $\Lambda$ of $\fa$ is invariant under the adjoint action of $T_1$ and $T_2$ on $\fa$. For $G=\Osc^2_q(\mu)$, this action is given by $e^{\rho(T_i)}\in \grO(\fa)$, $i=1,2$. For $G=D_q(\mu)$,  it equals $e^{\rho(T_i)}\oplus\Id_{\fa_0}\in \grO(\fa\oplus\fa_0)$, see Eq.\,(\ref{Elhl}).  This implies $e^{m_i\rho(T_i)}=\Id$ for suitable $m_i\in\NN\setminus\{0\}$, $i=1,2$. Hence, for $i=1,2$,  $\mu_1(T_i),\dots,\mu_q(T_i)$ span a lattice in $\RR$, which is generated by $2\pi q_i$ for some $q_i\in\QQ\setminus\{0\}$. Define $\sigma_1,\sigma_2\in\fz=(\RR^2)^*$ by $\sigma_i(T_j)=2\pi \delta_{ij}$. Then 
$$ \textstyle{\mu_j=\frac 1{2\pi} \mu_j(T_1)\sigma_1+\frac 1{2\pi} \mu_j(T_2)\sigma_2\in \ZZ\cdot q_1\sigma_1+ \ZZ\cdot q_2\sigma_2},$$
thus $\mu_j$ is in the lattice spanned by $q_1\sigma_1$ and $q_2\sigma_2$ for $j=1,\dots,q$.
\qed

\small\noindent Ines Kath\\Institut f\"ur Mathematik und Informatik, Universit\"at Greifswald,
Walther-Rathenau-Str. 47, \\ D-17487 Greifswald, Germany. \\
\texttt{ines.kath@uni-greifswald.de}
\end{document}